\documentclass[a4paper,10pt]{article}
\usepackage{amssymb,amsmath,amsthm}
\usepackage{fullpage}
\usepackage{hyperref}
\usepackage{eucal}
\usepackage{color}
\usepackage{stackrel}
\usepackage{graphicx}
\usepackage[shortlabels]{enumitem}

\newcommand{\etc}{\emph{etc.}}

\newcommand{\Real}{\mathbb{R}}
\newcommand{\Com}{\mathbb{C}}

\newcommand{\Dom}{\mathsf{D}}
\newcommand{\Num}{\mathsf{Num}}

\newcommand{\dist}{\mathop{\mathrm{dist}}\nolimits}

\newcommand{\sii}{L^2}

\newcommand{\der}{\mathrm{d}}

\newtheorem{Theorem}{Theorem}
\newtheorem{Lemma}{Lemma}

\theoremstyle{definition}

\def\OMIT#1{}
\usepackage[normalem]{ulem}
\definecolor{DarkGreen}{rgb}{0,0.5,0.1} 

\newcommand\soutD{\bgroup\markoverwith
{\textcolor{DarkGreen}{\rule[.5ex]{2pt}{1pt}}}\ULon}
\newcommand\soutP{\bgroup\markoverwith
{\textcolor{blue}{\rule[.5ex]{2pt}{1pt}}}\ULon}
\newcommand{\Hm}[1]{\leavevmode{\marginpar{\tiny%
$\hbox to 0mm{\hspace*{-0.5mm}$\leftarrow$\hss}%
\vcenter{\vrule depth 0.1mm height 0.1mm width \the\marginparwidth}%
\hbox to
0mm{\hss$\rightarrow$\hspace*{-0.5mm}}$\\\relax\raggedright #1}}}

\begin{document}
%
\title{\textbf{\LARGE On pseudospectra of 
the relativistic rotated harmonic oscillator}}
\author{David Krej\v{c}i\v{r}{\'\i}k\,$^a$ 
\ and \ Petr Siegl\,$^b$}
\date{\small 
\emph{
\begin{quote}
\begin{itemize}
\item[$a)$] 
Department of Mathematics, Faculty of Nuclear Sciences and 
Physical Engineering, Czech Technical University in Prague, 
Trojanova 13, 12000 Prague 2, Czechia;
david.krejcirik@fjfi.cvut.cz
\item[$b)$]
Institute of Applied Mathematics, 
Graz University of Technology, Austria;
siegl@tugraz.at
\end{itemize}
\end{quote}
}
\smallskip
11 November 2026}
\maketitle
 
%
\begin{abstract}
\noindent
We establish an upper bound to the norm of the resolvent 
of a relativistic version of the non-self-adjoint operator 
obtained by a dilation analytic
transformation of the quantum harmonic oscillator.
The estimate implies a decay for large complex energies
along rays with angles larger than that of the dilation.  
In this way we positively answer a question left open 
in~\cite{BKP}. 
 
%
%
%
\end{abstract}

\section{Introduction}
Almost any system in the universe 
that is disturbed slightly from a stable position 
behaves like a harmonic oscillator.
Perhaps everyone's earliest related experience is the swing.
Then we learn at school about classical examples
of vibrating spring, pendulum, 
inductor-capacitor electrical circuit, \etc.
In non-relativistic quantum mechanics, 
the harmonic oscillations are mathematically modelled by
the differential operator   
$
  S_0 := - \partial_x^2 + x^2 
$ 
in the space of square-integrable functions $\sii(\Real)$. 
The present paper is inspired by this fundamental physical system
enriched by non-self-adjointness and relativistic effects.
Both have been justified, around the turn of the millennium, 
by the advent of artificial materials like graphene 
and unconventional representations of quantum-mechanical observables.

To be more specific, we build upon 
the celebrated Davies' oscillator~\cite{Davies_1999a}
\begin{equation}\label{Davies}
  S_\theta := - e^{-i\theta} \partial_x^2 + e^{i\theta} x^2 
  \,,
  \qquad  
  \Dom({S_\theta}) := 
  H^2(\Real) \cap \sii(\Real,x^4 \, \der x)
  \,,
\end{equation}
in the Hilbert space $\sii(\Real)$,
where $\theta \in (-\frac{\pi}{2},\frac{\pi}{2})$.
The case $\theta=0$ represents 
the aforementioned Hamiltonian of the (self-adjoint)
harmonic oscillator in non-relativistic quantum mechanics.
The other values can be understood as obtained after a formal dilation transformation
$x \mapsto e^{i\theta/2} x$.
While the spectrum is preserved, 
$
  \sigma(S_\theta)
  = \sigma(S_0) = \{2n+1\}_{n=0}^\infty
$,
there is a drastic transition in finer (often called pseudospectral) 
properties of~$S_\theta$ whenever $\theta \not= 0$.   
In particular,
\begin{equation}\label{Davies.wild}
  \lim_{r \to +\infty}
  \big\|\big(S_\theta-r e^{i\vartheta}\big)^{-1}\big\|  
  \ = 
  \begin{cases}
  \infty
  & \mbox{if} \quad
  0 < |\vartheta| < |\theta|
  \,,
  \\
  0 
  & \mbox{if} \quad
  |\theta| < |\vartheta| \leq \pi
  \,,
  \end{cases}
\end{equation}
where~$\vartheta$ denotes the (real) argument of 
the spectral parameter (called energy) $r e^{i\vartheta}$.
Among the many important investigations of~$S_\theta$, 
let us point out  
\cite{Boulton_2002,Davies-Kuijlaars_2004,Pravda-Starov_2006,
Hitrik-Sjostrand-Viola_2013,Viola_2013,Henry,
KS4,Arnal-Siegl_2023} 
and the survey~\cite{KSTV}.

In~\cite{BKP}, a relativistic version of~$S_\theta$ 
was introduced in terms of the Dirac operator 
\begin{equation}\label{operator}
  H_\theta := - i \alpha_1 e^{-i\theta/2} \partial_x
  - \alpha_2 e^{i\theta/2} x
  + m \alpha_3
  \,,
  \qquad  
  \Dom({H_\theta}) := 
  [H^1(\Real) \cap \sii(\Real,x^2 \, \der x)]^4
  \,,
\end{equation}
in the Hilbert space $\sii(\Real)^4$.
Here~$m$ is a non-negative parameter (mass) 
and $\alpha_1,\alpha_2,\alpha_3$ denote
the standard $4 \times 4$ Dirac matrices
\begin{equation}\label{eq:alphamatrices}
  \alpha_k := 
  \begin{pmatrix}
    0 & \sigma_k \\
    \sigma_k & 0
  \end{pmatrix} 
  \,,
\end{equation}
with $k \in \{1,2,3\}$,
where
\begin{equation}\label{eq:sigmamatrices}
  \sigma_1 := 
  \begin{pmatrix}
    0 & 1 \\
    1 & 0
  \end{pmatrix} 
  \,, \qquad
  \sigma_2 := 
  \begin{pmatrix}
    0 & -i \\
    i & 0
  \end{pmatrix} 
  \,, \qquad
  \sigma_3 := 
  \begin{pmatrix}
    1 & 0 \\
    0 & -1
  \end{pmatrix} 
\end{equation}
are the Pauli matrices.
Using the supersymmetric relationship
\begin{equation}\label{eq:square}
  H_\theta^2 = (S_\theta + m^2) I_{\Com^4} + i \alpha_1\alpha_2 
  \,,
\end{equation}
where the matrix
\begin{equation}\label{eq:involutionmatrix}
  i \alpha_1\alpha_2  = 
  \begin{pmatrix}
    -\sigma_3 & 0  \\
    0 & -\sigma_3
  \end{pmatrix} 
\end{equation}
is diagonal,
it is possible to deduce that
\begin{equation}\label{eq:spectrum_Dirac}
  \sigma(H_\theta) = \sigma(H_0)
  = \left\{-\sqrt{2n+m^2},\sqrt{2n+m^2}\right\}_{n=0}^\infty
  \,.
\end{equation}
That is, up to the mass parameter, the spectrum of the non-relativistic~$S_\theta$
is a square of the spectrum of the relativistic~$H_\theta$, 
which is physically expected.
Moreover, the authors of~\cite{BKP} proved 
the related pseudospectral result 
\begin{equation}\label{wild} 
  \lim_{r \to \pm\infty}
  \big\|\big(H_\theta \mp m -r e^{i\vartheta}\big)^{-1}\big\|  
  = \infty
  \qquad \mbox{if} \qquad
  0 < |\vartheta| <  \mbox{$\frac{|\theta|}{2}$}
  \,.
\end{equation}
The question whether the half-angle is optimal was left open in~\cite{BKP}.
In the present paper, we fill in the gap by showing that 
the pseudospectral estimate is indeed sharp, as a consequence 
of the following result. 

\begin{Theorem}\label{Thm.main}
If $\frac{|\theta|}{2}  < |\vartheta| \leq  \frac{\pi}{2}$, then
$$
  \lim_{r \to \pm\infty}
  \big\|\big(H_\theta \mp m -r e^{i\vartheta}\big)^{-1}\big\|  
  = 0
  \,.
$$
\end{Theorem}

The rest of the paper consists of a proof of Theorem~\ref{Thm.main}.

\section{The proof}
Throughout the paper, we assume
$\theta \in (-\frac{\pi}{2},\frac{\pi}{2})$.
The proof is based on the following inequality 
for the non-relativistic oscillator.
\begin{Lemma}\label{Lem.Erice}
Let $z \in \Com \setminus \sigma(S_0)$ be non-zero.
Then
$$
  \cos(\theta)
  \left(
  \|\partial_x (S_\theta-z)^{-1}\|^2
  + \|x (S_\theta-z)^{-1}\|^2
  \right)
  \leq \frac{1}{|z|} 
  + 3 \, |z| \, \|(S_\theta-z)^{-1}\|^2
  \,.
$$
\end{Lemma}
\begin{proof}
For every $\psi \in \Dom(S_\theta) = \Dom(S_0)$,
one has 
$$
\begin{aligned}
  \|\psi\| \, \|(S_\theta-z)\psi\|
  \geq
  |(\psi,(S_\theta-z)\psi)|
  &= \left|
  e^{-i\theta} \|\partial_x\psi\|^2 + e^{i\theta} \|x\psi\|^2 
  - z \, \|\psi\|^2
  \right|
  \\
  &\geq   
  \left|
  e^{-i\theta} \|\partial_x\psi\|^2 + e^{i\theta} \|x\psi\|^2
  \right|
  - |z| \, \|\psi\|^2
  \\
  &\geq 
  \frac{|\cos(\theta)|}{2} \,  (\|\partial_x\psi\|^2 + \|x\psi\|^2)
  - |z| \, \|\psi\|^2
  \,.
\end{aligned}  
$$  
Consequently,
$$
  \cos(\theta) \,  (\|\partial_x\psi\|^2 + \|x\psi\|^2)
  \leq 2 \, |z| \, \|\psi\|^2
  + \delta \, \|\psi\|^2 + \delta^{-1} \, \|(S_\theta-z)\psi\|^2
$$
with every positive~$\delta$. 
Choosing $\delta:=|z|$ and $\psi := (S_\theta-z)^{-1}\phi$ 
with an arbitrary $\phi \in \sii(\Real)$,
we arrive at the desired claim.
\end{proof}

The following lemma particularly quantifies the decay of~\eqref{Davies.wild}.
Here and in the sequel, $\arg(z) \in (-\pi,\pi]$
denotes the prncipal branch of a complex number~$z$.
\begin{Lemma}\label{Lem.help}
Let $z \in \Com$ be sych that $|\theta| < |\arg(z)| \leq \frac{\pi}{2}$.
Then
$$
\begin{aligned}
  \|(S_\theta-z)^{-1}\| 
  &\leq \frac{1}{|z|} \frac{1}{\sin|\arg(z)-\theta|}
  \,, 
  \\
  \|\partial_x (S_\theta-z)^{-1}\|^2
  + \|x (S_\theta-z)^{-1}\|^2
  &\leq \frac{1}{|z|} 
  \frac{1}{\cos(\theta)}
  \left(
  1
  + \frac{3}{\sin^2|\arg(z)-\theta|}
  \right).
\end{aligned}  
$$
\end{Lemma}
\begin{proof}
The first inequality follows from the general property
$
  \|(S_\theta-z)^{-1}\| \leq 1/\dist\big(z,\overline{\Num(S_\theta)}\big) 
$,
where 
$
  \Num(S_\theta) 
  := \{(\psi,S_\theta \psi) : \psi \in \Dom(S_\theta)
  \ \land \ \|\psi\|=1 \}
$
is the numerical range of~$S_\theta$ and the fact that
$
  \overline{\Num(S_\theta)}
  \subset \{z \in \Com : |\arg(z)| \leq |\theta| \}
$  
in our case.
The second inequality then follows by Lemma~\ref{Lem.Erice}.
\end{proof}

The second ingredient in the proof is the following formula
for the resolvent of the Dirac operator.
\begin{Lemma}
Let $\lambda \in \Com \setminus \sigma(H_0)$. Then
$$
  (H_\theta -\lambda)^{-1} =
  \begin{pmatrix}
    \lambda R_- & 0 & m R_- & D R_+ \\
    0 & \lambda R_+ & D R_- & - m R_+ \\
    m R_- & D R_+ & \lambda R_- & 0 \\
    D R_- & -m R_+ & 0 & \lambda R_+ 
  \end{pmatrix}
  ,
$$
where
$$
  R_\pm := (S_\theta + m^2 - \lambda^2 \pm 1)^{-1}
  \qquad \mbox{and} \qquad
  D := -i e^{-i\theta/2} \partial_x + i e^{i\theta/2} x
  \,.
$$
\end{Lemma}
\begin{proof}
Using~\eqref{eq:square}, it follows that
$
  (H_\theta - \lambda)^{-1} 
  = (H_\theta + \lambda) 
  \left(
  (S_\theta + m^2-\lambda^2) I_{\Com^4} + i \alpha_1\alpha_2
  \right)^{-1}
$,
which is equivalent to the desired formula
when writing down explicitly using the Dirac matrices. 
\end{proof}
\begin{proof}[Proof of Theorem~\ref{Thm.main}]
The rest of our proof of Theorem~\ref{Thm.main}
is based on using the matrix max norm
to reduce the estimate of the norm of $(H_\theta -\lambda)^{-1}$ in $\sii(\Real)^4$
to routine bounds of the norm of the individual entries in $\sii(\Real)$
with the help of Lemma~\ref{Lem.help}.
For instance, writing $\lambda := m + r e^{i\vartheta}$
with $\frac{|\theta|}{2} < |\vartheta| \leq \frac{\pi}{2}$
and~$r$ large positive, one can estimate
$$
\begin{aligned}
  \|DR_-\|^2 
  &\leq 
  2 \, \|\partial_x (S_\theta - r^2 e^{i2\vartheta} - 2mr e^{i\vartheta} - 1)^{-1}\|^2
  + 2 \, \|x (S_\theta - r^2 e^{i2\vartheta} - 2mr e^{i\vartheta} - 1)^{-1}\|^2
  \\
  &\leq \frac{2}{|z|} 
  \frac{1}{\cos(\theta)}
  \left(
  1
  + \frac{3}{\sin^2|\arg(z)-\theta|}
  \right)
\end{aligned}  
$$
with $z := r^2 e^{i2\vartheta} + 2mr e^{i\vartheta} + 1$.
It is straightfoward to check that $|z| \sim r^2$ 
and $\arg(z) \sim 2\vartheta$ as $r \to +\infty$,
so Lemma~\ref{Lem.help} indeed applies
and $\|DR_-\| \to 0$ as $r \to \infty$.
At the same time,
$$
\begin{aligned}
  \|\lambda R_-\| 
  &= |m + r e^{i\vartheta}| \,
  \|(S_\theta - r^2 e^{i2\vartheta} - 2mr e^{i\vartheta} - 1)^{-1}\|
  \\
  &\leq |m + r e^{i\vartheta}| \,
  \frac{1}{|z|} \frac{1}{\sin|\arg(z)-\theta|}
  \to 0
\end{aligned} 
$$
as $r \to 0$.
The decay of the other terms can be checked similarly.
The complementary choice $\lambda := -m + r e^{i\vartheta}$ with $r \to -\infty$
in Theorem~\ref{Thm.main} follows analogously, too.
\end{proof}

\subsection*{Acknowledgments}
This article is based upon work from COST Action CA24122
(mSPACE -- multiscale Stochastics, Patterns, 
and Analysis of Combinatorial Environments)
supported by the European Cooperation in Science and Technology.
The first author was supported  
by the grant no.\ 26-21940S of the Czech Science Foundation.

\subsection*{AI usage}
None.
 
%
%

\providecommand{\bysame}{\leavevmode\hbox to3em{\hrulefill}\thinspace}
\providecommand{\MR}{\relax\ifhmode\unskip\space\fi MR }
\providecommand{\MRhref}[2]{%
  \href{http://www.ams.org/mathscinet-getitem?mr=#1}{#2}
}
\providecommand{\href}[2]{#2}

\end{document}